\documentclass[12pt]{amsart}

\newtheorem{theorem}{Theorem}[section]

\newtheorem{lemma}[theorem]{Lemma}
\theoremstyle{proposition}

\newtheorem{corollary}[theorem]{Corollary}
\usepackage[top=1.5in, bottom=1.5in, left=1.5in, right=1.5in]{geometry}
\newtheorem{definition}[theorem]{Definition}
\newtheorem{example}[theorem]{Example}

\theoremstyle{remark}

\newtheorem{remark}[theorem]{Remark}

\numberwithin{equation}{section}

\usepackage{times}
\usepackage{enumerate}
\usepackage{graphicx,adjustbox}
\usepackage{hyperref}
\usepackage{amsmath,amssymb}
\usepackage{amscd}
\usepackage{graphicx}
\usepackage[all]{xy}
\usepackage{mathtools}

\title[Associated Graded ring of Semigroup Algebras]
{On the Associated Graded ring of Semigroup Algebras}
\author{
Joydip Saha
\and
Indranath Sengupta
\and
Pranjal Srivastava
}
\date{}
\address{\small \rm  Stat Math Unit, Indian Statistical Institute, Kolkata, West-Bengal,700108, INDIA.} 
\email{saha.joydip56@gmail.com}
\thanks{The first author thanks NBHM, Government of India for post-doc fellowship at ISI kolkata}

\address{\small \rm  Discipline of Mathematics, IIT Gandhinagar, Palaj, Gandhinagar, 
Gujarat 382355, INDIA.}
\email{indranathsg@iitgn.ac.in}
\thanks{The second author is the corresponding author.}

\address{\small \rm  Discipline of Mathematics, IIT Gandhinagar, Palaj, Gandhinagar, 
Gujarat 382355, INDIA.}
\email{pranjal.srivastava@iitgn.ac.in}
\date{}



\begin{document}

\begin{abstract}
In this paper we give a necessary and sufficient condition for 
the Cohen-Macaulayness of the associated graded ring of a 
simplicial affine semigroups using Gr\"{o}bner basis. We generalize 
the concept of homogeneous numerical semigroup for the simplicial 
affine semigroup and show that the Betti numbers of the corresponding 
semigroup ring matches with the Betti numbers of the associated graded 
ring. We also define the nice extension for simplicial affine semigroups, 
motivated by the notion of a nice extension of the numerical semigroups.
\end{abstract}

\maketitle

\section{Introduction}

Let $S$ be an affine semigroup, fully embedded in $\mathbb{N}^{d}$. 
The semigroup algebra $k[S]$ over a field $k$ is generated by the 
monomials $x^{a}$, where $a \in S$, with maximal ideal 
$\mathfrak{m}=(x^{a_{1}},\dots,x^{a_{d+r}})$. Suppose that 
$S$ is a simplicial affine semigroup minimally generated by 
$\{a_{1}.\dots,a_{d},a_{d+1},\dots,a_{d+r}\}$, with the set 
of extremal rays $E=\{a_{1},\dots,a_{d}\}$. Many authors 
have studied the properties of the affine semigroup ring $k[S]$ from 
the properties of the affine semigroup $S$; see 
\cite{initial},\cite{Affine}.
\medskip

Let $I(S)$ denote the defining ideal of $k[S]$, which is the 
kernel of the $k-$algebra homomorphism  
$\phi:A=k[z_{1},\dots,z_{d+r}] \rightarrow k[x^{a_{1}},\dots,x^{a_{d+r}}]$, 
such that $\phi(z_{i})=x^{a_{i}}$, $i=1,\dots,d+r$. Let 
us write $k[S]\cong A/I(S)$. The defining ideal $I(S)$ 
is a binomial prime ideal (\cite{Herzog}, Proposition 1.4). 
The associated graded ring 
$\mathrm{gr}_{\mathfrak{m}}(k[S])=\oplus_{i=0}^{\infty}\mathfrak{m}^{i}/\mathfrak{m}^{i+1}$ is 
isomorphic to $\frac{k[z_{1},\dots,z_{d+r}]}{I(S)^{*}}$ 
(see \cite{initial}, Example 4.6.3), where $I(S)^{*}$ is the 
homogeneous ideal generated by the initial forms $f^{*}$ of the elements 
$f\in I(S)$, and $f^{*}$ is the homogeneous summand of $f$ of the 
least degree.
\medskip

Arslan et al. \cite{Gluing} have given the Gr\"{o}bner basis criterion for 
Cohen-Macaulayness of the associated graded ring of numerical 
semigroup rings and have used that to produce many examples which support 
Rossi's Conjecture. Herzog-Stamate have studied the Cohen-Macaulayness 
of the projective closure of monomial curves using Gr\"{o}bner basis 
(see \cite{Proj}). The projective closure of a numerical semigroup 
happens to be an affine semigroup in $\mathbb{N}^{2}$ and the 
associated graded ring is isomorphic to the semigroup 
ring. The homogeneity of the ideal $I(S)$ is the main property that 
is present in case of a projective closure of a numerical semigroup, 
which is not the case in general for an arbitrary affine semigroup. 
In this paper, we extend the result of Herzog-Stamate 
to a simplicial affine semigroup ring $k[S]$ and its associated graded ring 
$\mathrm{gr}_{\mathfrak{m}}(k[S])$. We prove, in Theorem \ref{CM}, that if 
$G=\{f_{1},\dots,f_{r}\}$ is the minimal Gr\"{o}bner basis of the 
defining ideal $I(S)$, with respect to the negative degree reverse 
lexicographic ordering induced by $z_{d+r} > \dots > z_{d} > \dots >z_{1}$, 
then $\mathrm{gr}_{\mathfrak{m}}(k[S])$ is Cohen-Macaulay if and only if 
$z_{j}$ does not divide the $\mbox{LM}(f_{i})$, for every 
$1\leq j \leq d, 1 \leq i \leq r $. 
\medskip

Numerical semigroups are important in the study of curve singularities 
and these are natural classes of simplicial affine 
semigroups. We have tried to extend some important 
results on numerical semigroups to the context of affine semigroups. Jafari-Zarzuela \cite{Homo type} 
have defined the noition of a homogeneous numerical 
semigroup $S=\langle n_{1},\dots, n_{r}\rangle$ as the one 
which has the property that every element of $\mathrm{AP}(S,n_{1})$ 
has a maximal expression, 
and have studied the Betti numbers of their associated graded rings. 
We show that the notion of homogeneity can be generalized for 
affine semigroups if we use 
the Ap\'{e}ry set of $S$ with respect to the set of extremal 
rays $E$ and we prove the following: If $\mathrm{gr}_{\mathfrak{m}}(k[S])$ of a homogeneous 
simplicial affine semigroup $S$ is Cohen-Macaulay, then the Betti numbers 
of $k[S]$ coincide with the Betti numbers of $\mathrm{gr}_{\mathfrak{m}}(k[S])$, i.e., 
$\beta_{i}(\mathrm{gr}_{\mathfrak{m}}(k[S])) = \beta_{i}(k[S])$. This observation 
is in the spirit of some of the earlier results proved for affine semigroups $S$, viz., 
the study of Cohen-Macaulayness of 
$k[S]$ and its associated graded ring $\mathrm{gr}_{\mathfrak{m}}(k[S])$ by  
Jafari et al.\cite{Reduction} and the observation due to Herzog et.al. \cite{Betti}, that, 
$\beta_{i}(\mathrm{gr}_{\mathfrak{m}}(k[S])) \geq \beta_{i}(k[S])$, for all $i \geq 1$.
\medskip

In \cite{Nice}, Feza Arslan defined the nice extension of numerical 
semigroups. We generalize this notion of nice extension to 
simplicial affine semigroups 
and prove that the  associated graded ring of the nice extension is always 
Cohen-Macaulay if the original semigroup is Cohen-Macaulay.  
\medskip

Let us discuss how the sections are divided in this article. In section 3, 
we give a Gr\"{o}bner basis criterion for the Cohen-Macaulayness 
of the associated graded ring, using minimal reduction ideal of 
the maximal ideal 
$\mathfrak{m}$. Section 4 is devoted to the study of homogeneous simplicial 
affine semigroups. We also give some examples of homogeneous simplicial affine 
semigroups and prove that if the associated graded ring $\mathrm{gr}_{\mathfrak{m}}(k[S])$ 
of the semigroup ring $k[S]$ associated to the homogeneous simplicial affine 
semigroup $S$ is Cohen-Macaulay then 
$\beta_{i}(\mathrm{gr}_{\mathfrak{m}}(k[S])) = \beta_{i}(k[S])$ (Theorem \ref{Betti}). In section 5, 
we explain the nice extension of simplicial affine semigroups and prove some of its 
properties like complete intersection and Cohen-Macaulayness. We also show that every 
semigroup ring associated with a simplicial affine semigroup, obtained by a sequence of nice 
extension of a simplicial affine semigroup, is a complete intersection (Theorem \ref{CI}).

\section{Preliminaries}
Let $S$ be an affine semigroup, fully embedded in $\mathbb{N}^{d}$. 

\begin{definition}{\rm 
The \textit{rational polyhedral cone} generated by $S$ is defined as
\[
 \mathrm{cone}(S)=\big\{\sum_{i=1}^{n}r_{i}a_{i}: r_{i} \in \mathbb{R}_{\geq 0}, \,i=1,\dots,d+r\big\}.
 \]
 The \textit{dimension} of $S$ is defined as the dimension of the subspace generated by $\mathrm{cone}(S)$.
 }
\end{definition} 
   
The $\mathrm{cone}(S)$ is the intersection of finitely many closed linear half-spaces in $\mathbb{R}^{d}$, 
each of whose bounding hyperplanes contains the origin. These half-spaces are called \textit{support hyperplanes}. 

\begin{definition} {\rm 
Suppose $S$ is an affine semigroup, fully embedded in $\mathbb{N}^{d}$. 
If $d=2$, the support hyperplanes are one-dimensional vector spaces, which are called the 
\textit{extremal rays} of $\mathrm{cone}(S)$. When $d >2$, intersection of any two
adjacent support hyperplanes is a one-dimensional vector space, called an extremal ray 
of $\mathrm{cone}(S)$. An element of $S$ is called an extremal ray of $S$ 
if it is the smallest non-zero vector of $S$ in an extremal ray of $\mathrm{cone}(S)$. 
}
\end{definition}

\begin{definition}\label{Extremal}{\rm 
An affine semigroup $S$, fully embedded in $\mathbb{N}^{d}$, is said to be \textit{simplicial}   
if the $\mathrm{cone}(S)$ has atleast $d$ extremal rays, i.e., if there exist $d$ elements say 
$\{a_{1},\dots,a_{d}\} \subset \{a_{1}.\dots,a_{d},a_{d+1},\dots,a_{d+r}\}$, such that they 
are linearly independent over $\mathbb{Q}$ and $S \subset \sum\limits_{i=1}^{d}\mathbb{Q}_{\geq 0}a_{i}$.
}
\end{definition}

In this paper, $S$ always denotes a simplicial affine semigroup minimally generated by 
$\{a_{1}.\dots,a_{d},a_{d+1},\dots,a_{d+r}\}$, with the set of extremal rays $E=\{a_{1}.\dots,a_{d}\}$. 
The semigroup ring defined by $S$ is written as $k[S]=k[x^{a_{1}},\dots,x^{a_{d+r}}]$.

\begin{remark}
We note that if $f=z^{p}-z^{q} \in I(S)$, where $p$ and $q$ are $d$-tuples of non-negative integers, the 
set $\{z_{j} \mid p_{j}+q_{j}\neq 0\}$ is the support of $f$ denoted by $\mathrm{supp}(f)$.
\end{remark}

\begin{definition}
{\rm 
A subset $H \subseteq S$ is called an ideal of $S$, if $H+S \subseteq H$. 
}
\end{definition}
Suppose $H_{1}$ and $H_{2}$ be two ideals of $S$. We define 
$H_{1}+H_{2}=\{h_{1}+h_{2} \mid h_{1}\in H_{1}, h_{2}\in H_{2}\}$. 
For a positive integer $n$ and an ideal $H$, we define $nH$ as $H+(n-1)H$ and $2H=H+H$. 
Let $M=S \setminus \{0\}$ be the maximal ideal of $S$. 

\begin{definition}\label{maximal}{\rm 
The maximum integer $n$, such that $s \in nM \setminus (n+1)M$, is called the \textit{order} of 
$s$, written as $n=\mathrm{ord}_{S}(s)$. If $s=\sum\limits_{i=1}^{d+r}r_{i}a_{i}$ for 
some non-negative integers 
$r_{i}$, such that $\sum\limits_{i=1}^{d+r}r_{i}=n = \mathrm{ord}_{S}(s)$, it is called a 
\textit{maximal expression} of $s$ and $(r_{1},\dots,r_{d+r})$ is called a \emph{maximal factorization} of $s$. 
}
\end{definition}

\begin{definition}\emph{
The Ap\'{e}ry set of $S$ with respect to an element $b \in S$ is defined as $\{a \in S: a-b \notin S\}$. Let $E=\{a_{1},\dots,a_{d}\}$ be a set of extremal rays of $S$, then the Ap\'{e}ry set of $S$ with respect to the set $E$ is
\[
\mathrm{AP}(S,E)=\{a \in S \mid a-a_{i} \notin S, \, \forall i=1,\dots,d\}=\cap_{i=1}^{d}\mathrm{AP}(S,a_{i})
\]}
\end{definition}

\begin{definition}{\rm
Let $A$ be a graded noetherian ring,  $I$ an ideal in $A$. Let $I^{*}$ be the ideal in $A$, 
generated by the element $f^{*}$, where $f\in I$ and $f^{*}$ is the homogeneous summand of $f$ 
with the least total degree. A set $\{f_{1},\dots,f_{t}\}\subseteq I$ is called a 
\textit{standard basis} for $I$ if $I^{*}$ is generated by $\{f_{1}^{*},\dots,f_{t}^{*}\}$. 
}
\end{definition}

\begin{definition}{\rm
Let $(B,\mathcal{F})$ be a filtered, Noetherian ring. A sequence $g = g_{1},\dots,g_{n}$ in $B$ is called 
\textit{super regular} if the sequence of initial forms $g^{*} = g_{1}^{*},\dots,g_{n}^{*}$ is 
regular in $\mathrm{gr}_{\mathcal{F}}(B)$.
}
\end{definition}

\section{Gr\"{o}bner Basis Criterion of Cohen-Macaulay of associated graded ring}
We prove a condition for the Cohen-Macaulayness of the associated graded ring $\mathrm{gr}_{\mathfrak{m}}(k[S])$, 
where $\mathfrak{m}=(x^{a_{1}},\dots,x^{a_{d+r}})$. The 
condition that we establish involves a Gr\"{o}bner basis of $I(S)$ and hence it is computational 
in nature.  Under some mild conditions on the Gr\"{o}bner basis (see Theorem \ref{criterion}), 
we also prove that the Betti sequence of the $\mathrm{gr}_{\mathfrak{m}}(k[S])$ is exactly the 
same as the Betti sequence of the semigroup ring $k[S]$. Let us discuss some lemmas first.

\begin{lemma}\label{gbasis}
We consider the following map $$\pi_{i}:k[z_{1},\dots,z_{d},\dots,z_{d+r}] \rightarrow k[\hat{z_{1}},\dots,\hat{z_{i}},z_{i+1},\dots,z_{d+r}],$$ such that $\pi_{i}(z_{j})=0, 1 \leq j \leq i$ and $\pi_{i}(z_{j})=z_{j}, i+1 \leq j \leq d+r$. 
Let
$G=\{f_{1},\dots,f_{t}\}$ be a Gr\"{o}bner basis of the defining ideal $I(S)$, with respect to the 
negative degree reverse lexicographic ordering induced by $z_{d+r} > \dots > z_{d} > \dots >z_{1}$ 
on $k[z_{1},\dots,z_{d},\dots,z_{d+r}]$. 
If for every $1 \leq j \leq i$, $z_{j}$ does not divide the leading monomial of any element of $G$, 
then $\pi_{i}(G)=\{\pi_{i}(f_{1}),\dots,\pi_{i}(f_{t})\}$ is a Gr\"{o}bner basis of $\pi_{i}(I(S))$, 
with respect to the negative degree reverse lexicographic ordering induced by 
$z_{d+r} > \dots > z_{d} > \dots >z_{i+1}$ on $k[\hat{z_{1}},\dots,\hat{z_{i}},z_{i+1},\dots,z_{d+r}]$.
\end{lemma}

\proof For all $ 1 \leq j \leq i$, $z_{j}$ do not divide leading monomial of any element of $G$. Therefore, 
we have $\pi_{i}(\mathrm{LM}(f_{l}))=\mathrm{LM}(\pi_{i}(f_{l}))$, where $f_{l}\in G$. Let 
$\pi_{i}(f) \in \pi_{i}(I(S))$, for some $f \in I(S)$ and $f=\mathrm{LM}(f)+g$, for $g \in A$. 
If $z_{j}$ divides $\mathrm{LM}(f)$, for some $j \in \{1,\dots,i\}$, then due to the negative 
degree reverse lexicographic ordering induced by $z_{d+r} > \dots > z_{d} > \dots >z_{1}$, either 
$z_{j}$ divides $g$ or $z_{l}$ divides $g$ for some $l <j$. Hence $\pi_{i}(f)=0$. Therefore, 
$$\mathrm{LM}(\pi_{i}(f))=0=\pi_{i}(\mathrm{LM}(f)) \in \pi_{i}(\langle f_{1},\dots,f_{t}\rangle)=\langle \pi_{i}(f_{1}),\dots,\pi_{i}(f_{t})\rangle.$$

If $z_{j} \nmid \mathrm{LM}(f))$, for any $j=1, \dots, i$, then 
$\pi_{i}(\mathrm{LM}(f))=\mathrm{LM}(f)$ and $\pi_{i}(f)=\mathrm{LM}(f)+\pi_{i}(g)$, f
or some $g \in A$. Therefore, 
$$\mathrm{LM}(\pi_{i}(f))=\pi_{i}(\mathrm{LM}(f)) \in \pi_{i}(\langle f_{1},\dots,f_{t}\rangle)=\langle \pi_{i}(f_{1}),\dots,\pi_{i}(f_{t})\rangle.$$ 
Hence, $\pi_{i}(G)=\{\pi_{i}(f_{1}),\dots,\pi_{i}(f_{t})\} $ is a Gr\"{o}bner basis of $\pi_{i}(I(S))$, 
with respect to the negative degree reverse lexicographic ordering induced by $z_{d+r} > \dots > z_{d} > \dots >z_{i+1}$. \qed

\begin{definition}{\rm
Let $B$ be a Noetherian ring, $I$ a proper ideal and $M$ a finite $B$-module. An ideal 
$J \subset I$ is called a \textit{reduction ideal} of $I$, with
respect to $M$, if $JI^{n}M =I^{n+1}M$ for some (or equivalently all) sufficiently large $n$. 
}
\end{definition}

\begin{lemma}\label{Cond}
Let $(x^{a_{1}},\dots,x^{a_{d}})$ be a reduction ideal of $\mathfrak{m}$, then the following statements are equivalent:
\begin{enumerate}[(a)]
\item $\mathrm{gr}_{\mathfrak{m}}(k[S])$ is a Cohen-Macaulay ring.
\item $(x^{a_{1}})^{*},\dots,(x^{a_{d}})^{*}$ provides a regular sequence in $\mathrm{gr}_{\mathfrak{m}}(k[S])$.
\item $R$ is Cohen-Macaulay and $(x^{a_{i}})^{*}$ is a non-zero divisor in 
$\mathrm{gr}_{\mathfrak{m}}(k[S]),\, \text{for} \,\, i=1,\dots,d$.
\end{enumerate}

\end{lemma}

\proof See Proposition 5.2 in \cite{Reduction}. \qed

\begin{remark}\label{image z}
We have the map $\phi:A=k[z_{1},\dots,z_{d+r}] \rightarrow k[x^{a_{1}},\dots,x^{a_{d+r}}]$, defined as 
$\phi(z_{i})=x^{a_{i}}$, for $i=1,\dots,d+r$. Therefore, $\phi$ is a surjective map and we have the 
induced surjective map $\mathrm{gr}(\phi):\mathrm{gr}_{\eta}(A)\rightarrow \mathrm{gr}_{\mathfrak{m}}(k[S])$, 
such that $\mathrm{gr}(\phi)(z_{i})=\phi(z_{i})+\mathfrak{m}^{2}=x^{a_{i}}+\mathfrak{m}^{2}=(x^{a_{i}})^{*}$, 
where $\eta = (z_{1}, \ldots , z_{d+r})$, $z_{i} \in \eta \setminus \eta^{2}$ and 
$\mathrm{gr}_{\mathfrak{m}}(k[S]) \cong A/I(S)^{*}$.  
\end{remark}

\begin{theorem}\label{criterion}
Let $S$ be a simplicial affine semigroup, fully embedded in $\mathbb{N}^{d}$, such that $k[S]$ is Cohen-Macaulay. 
Let $(x^{a_{1}},\dots,x^{a_{d}})$ be a reduction ideal of $\mathfrak{m} = (x^{a_{1}},\dots,x^{a_{d+r}})$. Let 
$G=\{f_{1},\dots,f_{t}\}$ be the minimal Gr\"{o}bner basis of the defining ideal $I(S)$, with respect to the 
negative degree reverse lexicographic ordering $z_{d+r} > \dots > z_{d} > \dots >z_{1}$. 
Then, $\mathrm{gr}_{\mathfrak{m}}(k[S])$ is Cohen-Macaulay if and only if for every 
$1\leq j \leq d$, $1 \leq i \leq t$, the indeterminate $z_{j}$ does not divide $\mathrm{LM}(f_{i})$.  
\end{theorem}
 
\proof We proceed by induction on $d$.

\noindent{\textbf{Case: $\mathbf{d=1}$.}}
Let $z_{1}$ divide $\mathrm{LM}(f_{i})$, for some $i$. Then $f_{i}^{*}=z_{1}m+\sum c_{i}m_{i}$, 
where $m_{i}$ are monomials such that $z_{1}$ divides each $m_{i}$. If $z_{1}$ fails to divide 
at least one $m_{i}$, the leading monomial will be that term where $z_{1}$ is not present. 
Therefore, $f_{i}^{*}=z_{1}g$, where $g$ is a homogeneous polynomial. If $g \in I(S)^{*}$, then 
$f^{*}=g$, for some $f \in I(S)$, and $\mathrm{LM}(f)=\mathrm{LM}(g)$. Since $G$ is a Gr\"{o}bner basis of $I(S)$, 
$\mathrm{LM}(f_{j})$ divides $\mathrm{LM}(f)=\mathrm{LM}(g)$, for some $f_{j} \in G$, 
and $\mathrm{LM}(g)$ divides $\mathrm{LM}(f_{i})$. These imply that 
$\mathrm{LM}(f_{j})$ divides $\mathrm{LM}(f_{i})$, which contradicts the minimality 
of $G$. Hence, $g \notin I(S)^{*}$. Therefore, $\bar{z_{1}}$ is a zero divisor in $A/I(S)^{*}$ 
and hence  $\mathrm{gr}(\phi)(z_{i})=(x^{a_{1}})^{*}$ is a zero divisor in $\mathrm{gr}_{\mathfrak{m}}(k[S])$, 
where $\mathrm{gr}(\phi)$ has been defined in Remark \ref{image z}. This proves that 
$\mathrm{gr}_{\mathfrak{m}}(k[S])$ is not Cohen-Macaulay. 

Conversely, if $A/I(S)^{*}$ 
is not Cohen-Macaulay, then $\bar{z_{1}}$ is a zero divisor in $A/I(S)^{*}$. Therefore, 
$z_{1}g \in I(S)^{*}$, where $g$ is a
monomial or a homogeneous polynomial with $g \notin I(S)^{*}$, and $\mathrm{LM}(f_{i})$ does 
not divide $\mathrm{LM}(g)$, for every $1 \leq i \leq t$. Since the ideal generated by the leading 
monomials of the elements
in $I(S)$ contains $z_{1}\mathrm{LM}(g)$, there exists $f_{i} \in G$ such that $\mathrm{LM}(f_{i}) = z_{1}m,$ where $m$ is a monomial that divides $\mathrm{LM}(g)$. 
\medskip

\noindent{\textbf{Case: $\mathbf{d\geq 2}$.}} We assume that the result holds for $1\leq d \leq j-1$. 
Now for induction step assume $z_{1},\dots,z_{j-1}$ does not divide $\mathrm{LM}(f_{i})$, 
but $z_{j} \mid \mathrm{LM}(f_{i})$ for some $i$. 

Consider the map $$\pi_{j-1}:k[z_{1},\dots,z_{d},\dots,z_{d+r}] \rightarrow \bar{A}:=k[z_{j},\dots,z_{d+r}],$$ 
such that $\pi_{j-1}(z_{l})=0$, $1 \leq l \leq j-1$ and $\pi_{j-1}(z_{l})=z_{l}$, $j \leq l \leq d+r$. 
Since $z_{1},\dots,z_{j-1}$ do not divide $\mathrm{LM}(f_{i})$, for all $i$, 
and $z_{j}$ divides $\mathrm{LM}(f_{i})$, for some $i$, then 
$f_{i}^{*}=z_{j}m+\sum c_{p}m_{p}$, where $m_{p}$ are monomials, 
and $\pi_{j-1}(f_{i}^{*})=z_{j}\pi_{j-1}(m)+\sum c_{p}\pi_{j-1}(m_{p})$. 
By Lemma \ref{gbasis}, $\pi_{j-1}(G)=\{\pi_{j-1}(f_{1}),\dots,\pi_{j-1}(f_{r})\}$ 
is a minimal Gr\"{o}bner basis of $\pi_{j-1}(I(S))$, with respect to the 
negative degree reverse lexicographic ordering induced by $z_{d+r} > \dots >z_{d+1} > \dots > z_{j}$. 
Moreover, $\mathrm{LM}(f^{*}) = \mathrm{LM}(f)$, due to 
the very choice of the monomial order defined above. Now, given that $z_{j} \mid \mathrm{LM}(f_{i})$, 
we must have that $z_{j} \mid \pi_{j-1}(m_{p})$ for each $p$, failing which, the leading 
monomial of $\pi_{j-1}(f_{i})$ comes from $\pi_{j-1}(m_{p})$, for some $p$, 
because of the choice of the monomial ordering, contradicting $z_{j} \mid \mathrm{LM}(f_{i})$. 
Therefore, $\pi_{j-1}(f_{i}^{*})=z_{j}\pi_{j-1}(g)$, for some homogeneous polynomial $g$. 
\medskip

If $\pi_{j-1}(g)\in \pi_{j-1}(I(S)^{*})$, then $\pi_{j-1}(f^{*})=\pi_{j-1}(g)$ for some $f \in I(S)$ and 
$\mathrm{LM}(\pi_{j-1}(f)) = \mathrm{LM}(\pi_{j-1}(f^{*})) = \mathrm{LM}(\pi_{j-1}(g))$. 
Since $\pi_{j-1}(G)$ is a Gr\"{o}bner basis of $\pi_{j-1}(I(S))$, we have  $\mathrm{LM}(\pi_{j-1}(f_{j}))$ 
divides $\mathrm{LM}(\pi_{j-1}(f))=\mathrm{LM}(\pi_{j-1}(g))$ for some $f_{j} \in G$.  We know that 
$\mathrm{LM}(\pi_{j-1}(g))$ divides $\mathrm{LM}(\pi_{j-1}(f_{i}))$. Therefore, 
$\mathrm{LM}(\pi_{j-1}(f_{j}))$ divides $\mathrm{LM}(\pi_{j-1}(f_{i}))$, which contradicts 
the minimality of $\pi_{j-1}(G)$. Hence, $\pi_{j-1}(g) \notin (\pi_{j-1}(I(S)^{*})$. 
Therefore, $z_{j}$ is a zero-divisor in $\bar{A}/\pi_{j-1}(I(S)^{*}) \cong A/(z_{1},\dots,z_{j-1},I(S)^{*})$ 
and hence $(x^{a_{i}})^{*}$ is a zero-divisor in $\mathrm{gr}_{\mathfrak{m}}(R)$. This proves that 
$\mathrm{gr}_{\mathfrak{m}}(R)$ is not Cohen-Macaulay. 
\medskip

Conversely, suppose $A/I(S)^{*}$ is not Cohen-Macaulay. By the induction hypothesis, 
we may assume that $z_{1},\dots,z_{j-1}$ form a regular sequence in $A/I(S)^{*}$ and 
$z_{j}$ is a zero-divisor in $ A/(z_{1},\dots,z_{j-1},I(S)^{*}) = \bar{A}/(\pi_{j-1}I(S)^{*})$. 
Moreover, $z_{1},\dots,z_{j-1}$ do not divide $\mathrm{LM}(f_{i})$, for any $i$. 
Therefore $\mathrm{LM}(\pi_{j-1}(f_{i})) = \mathrm{LM}(f_{i})$, for all $i$. 
From Proposition 15.13 in \cite{Eisenbud} we have  $z_{j}$ is a zero-divisor in $\bar{A}/\mathrm{LI}(\pi_{j-1}I(S)^{*})$, where $\mathrm{LI}(\pi_{j-1}I(S)^{*})$ denotes the 
leading ideal of $\pi_{j-1}(I(S)^{*})$.  Therefore, 
$z_{j}g \in \mathrm{LI}(\pi_{j-1}I(S)^{*})$, for some monomial $g \notin \mathrm{LI}(\pi_{j-1}I(S)^{*})$. 
This implies $ \mathrm{LM}(\pi_{j-1}(f_{i})) \nmid g$  for any $i$. 
But $z_{j}g\in \mathrm{LI}(\pi_{j-1}(I(S))) = \mathrm{LI}(\pi_{j-1}I(S)^{*})$, which implies that 
$z_{j}g=m'\mathrm{LM}(\pi_{j-1}(f_{i}))$, for some $f_{i}\in G$ and a monomial $m'$. Suppose $z_{j} \nmid \mathrm{LM}(\pi_{j-1}(f_{i}))$ for all $i$, then $z_{j}\mid m'$ and we get $ \mathrm{LM}(\pi_{j-1}(f_{i})) \mid g$, which is a contradiction. Hence $z_{j}\mid \mathrm{LM}(\pi_{j-1}(f_{i}))$. We know 
that $\mathrm{LM}(\pi_{j-1}(f_{i})) = \mathrm{LM}(f_{i})$, hence 
$z_{j} \mid \mathrm{LM}(f_{i})$ for all $i$.
 \qed
\medskip

Let $A$ be a filtered noetherian graded ring with homogeneous maximal ideal $\mathfrak{m}_{A}$ and suppose $B=A/xA$, where $x$ is not a zero-divisor on $A$. Let 
$\psi:A \rightarrow B$ be the canonical epimorphism.

\begin{lemma}\label{iso}
If $x$ is a super regular in $A$ then 
\[
\mathrm{gr}_{\mathfrak{m}_{A}}(A) \xrightarrow{(x)^{*}} \mathrm{gr}_{\mathfrak{m}_{A}}(A) \xrightarrow{\mathrm{gr}(\psi)} \mathrm{gr}_{\mathfrak{m}_{B}}(B) \rightarrow 0
\]
is exact. 
\end{lemma} 
\proof See Lemma a in \cite{Super}. \qed

\begin{lemma}\label{iso1}
Consider a map $$\pi_{d}:k[z_{1},\dots,z_{d},\dots,z_{d+r}] \rightarrow \bar{A}=k[z_{d+1},\dots,z_{d+r}]$$ such that $\pi_{d}(z_{j})=0, 1 \leq j \leq d$ and $\pi_{d}(z_{j})=z_{j}, d+1 \leq j \leq d+r$. If $z_{1},\dots,z_{d}$ is a super regular in $A/I(S)$ then

$$\mathrm{gr}_{\bar{\mathfrak{m}}}\big(\bar{A}/\pi_{d}(I(S)) \cong \frac{\mathrm{gr}_{\mathfrak{m}}(A/I(S))}{(z_{1},\dots,z_{d})\mathrm{gr}_{\mathfrak{m}}(A/I(S))},$$ 
where  $\bar{\mathfrak{m}}=\pi_{d}(\mathfrak{m})$.
\end{lemma}

\proof Consider the exact sequence 
\[
0 \xrightarrow{} (z_{1},\dots,z_{d})\frac{A}{I(S)}  \xrightarrow{\mathrm{ker}(\pi_{d})} \frac{A}{I(S)} \xrightarrow{\pi_{d} } \frac{\bar{A}}{\pi_{d}(I(S))} \rightarrow 0.
 \]
Since  $z_{1},\dots,z_{d}$ is super regular in $A/I(S)$, by Lemma \ref{iso}, we have an exact sequence
 \[
0 \xrightarrow{} (z_{1},\dots,z_{d})\mathrm{gr}_{\mathfrak{m}}\Big(\frac{A}{I(S)} \Big) \xrightarrow{\mathrm{gr}_{\mathfrak{m}}(\mathrm{ker}\pi_{d})} \mathrm{gr}_{\mathfrak{m}}\Big(\frac{A}{I(S)} \Big)\xrightarrow{\mathrm{gr}_{\mathfrak{m}}(\pi_{d}) } \mathrm{gr}_{\bar{\mathfrak{m}}}\Big(\frac{\bar{A}}{\pi_{d}(I(S))}\Big) \rightarrow 0
\].
Therefore,$$\mathrm{gr}_{\bar{\mathfrak{m}}}\big(\bar{A}/\pi_{d}(I(S)) \cong \frac{\mathrm{gr}_{\mathfrak{m}}(A/I(S))}{(z_{1},\dots,z_{d})\mathrm{gr}_{\mathfrak{m}}(A/I(S))}. \qed$$ 
\medskip
 
\begin{theorem}\label{CM}
Let $(x^{a_{1}},\dots,x^{a_{d}})$ be a reduction ideal of $\mathfrak{m}$. 
Suppose $k[S]$ and $\mathrm{gr}_{\mathfrak{m}}(k[S])$ are Cohen-Macaulay. 
Let $G=\{f_{1},\dots,f_{t},g_{1},\dots,g_{s}\}$ be a minimal Gr\"{o}bner basis of the defining ideal $I(S)$, with respect to the negative degree reverse 
lexicographic ordering induced by $z_{d+r} > \dots > z_{d} > \dots >z_{1}$. 
We assume that $f_{1}, \dots,f_{t}$ are homogeneous and $g_{1},\dots,g_{s}$ 
are non homogeneous, with respect to the standard gradation on the polynomial 
ring $k[z_{1}, \ldots, z_{d+r}]$. If there exists a $j$, $1\leq j \leq d$, such that 
$z_{j}$ belongs to the support of $g_{l}$, for every $1 \leq l \leq s $, 
then $$\beta_{i}(k[S])=\beta_{i}(\mathrm{gr}_{\mathfrak{m}}(k[S])) 
\quad \forall i \geq 1.$$
\end{theorem}

\proof 
Let $G=\{f_{1},\dots,f_{r},g_{1},\dots,g_{s}\}$ be a minimal Gr\"{o}bner basis of the defining ideal $I(S)$, with respect to the negative degree reverse 
lexicographic ordering induced by $z_{d+r} > \dots > z_{d} > \dots >z_{1}$. 
When $s=0$, $I(S)$ is homogeneous ideal and from Remark 2.1 (\cite{Reduction}),  
$k[S] \cong \mathrm{gr}_{\mathfrak{m}}(k[S])$. Hence, the 
result follows directly.
\medskip

When $s\geq 1$, we have $f_{1}, \dots, f_{t}$ are homogeneous, 
$g_{1},\dots g_{s}$ are non-homogeneous and $\mathrm{gr}_{\mathfrak{m}}(k[S])$ 
is Cohen-Macaulay, this implies that $z_{1},\dots,z_{d}$ do not divide the 
$\mathrm{LM}(f_{k})$ 
and $\mathrm{LM}(g_{l})$ for $k=1,\dots, t$, $l=1,\dots, s$. Moreover, 
$z_{j} \in \mathrm{supp}(\{g_{1},\dots,g_{s}\})$, for some $1 \leq j \leq d$. Therefore, $z_{j}$ divides a non-leading term of $g_{1},\dots,g_{s}$, for 
some  $1 \leq j \leq d$.
\medskip

We consider the map $$\pi_{d}:A=k[z_{1},\dots,z_{d},\dots,z_{d+r}] \rightarrow \bar{A}=k[z_{d+1},\dots,z_{d+r}]$$ such that $\pi_{d}(z_{j})=0$, 
$1 \leq j \leq d$ and $\pi_{d}(z_{j})=z_{j}, d+1 \leq j \leq d+r$. 
We note that $\pi_{d}(f_{1}),\dots,\pi_{d}(f_{r})$ are either monomials 
or homogeneous polynomials. Since $z_{j}$ divides a non-homogeneous 
term of $\{g_{1},\dots,g_{s}\}$ for some $1 \leq j \leq d+r$, we must have 
that $\pi_{d}(g_{1}),\dots,\pi_{d}(g_{s})$ are the leading monomials of $g_{1},\dots,g_{s}$ respectively. 
\medskip

Therefore $\{\pi_{d}(f_{1}),\dots,\pi_{d}(f_{r}),\pi_{d}(g_{1}),\dots,\pi_{d}(g_{s})\}$ generates the homogeneous ideal $\pi_{d}(I(S))$. Hence
\begin{align*}
\beta_{i}\big(\bar{A}/\pi_{d}(I(S))\big)=& \beta_{i}(\mathrm{gr}_{\bar{\mathfrak{m}}}\big(\bar{A}/\pi_{d}(I(S))\big) 
\end{align*}
where  $\bar{\mathfrak{m}}=\pi_{d}(\mathfrak{m})$. By 
Lemma \ref{iso1}, 
\[\mathrm{gr}_{\bar{\mathfrak{m}}}\big(\bar{A}/\pi_{d}(I(S)) \cong \frac{\mathrm{gr}_{\mathfrak{m}}(A/I(S))}{(z_{1},\dots,z_{d})\mathrm{gr}_{\mathfrak{m}}(A/I(S))},
\]  
 therefore $$\beta_{i}(\mathrm{gr}_{\bar{\mathfrak{m}}}\big(\bar{A}/\pi_{d}(I(S))\big) = \beta_{i}\bigg(\frac{\mathrm{gr}_{\mathfrak{m}}(A/I(S))}{(z_{1},\dots,z_{d})\mathrm{gr}_{\mathfrak{m}}(A/I(S))}\bigg).$$
$\mathrm{gr}_{\mathfrak{m}}(A/I(S))$ being Cohen-Macaulay, $z_{1},\dots,z_{d}$ form 
a regular sequence in $\mathrm{gr}_{\mathfrak{m}}(A/I(S))$, by Lemma \ref{Cond}. 
We know that the Betti numbers are preserved under going modulo regular 
elements, hence
$$ \beta_{i}\bigg(\frac{\mathrm{gr}_{\mathfrak{m}}(A/I(S))}{((z_{1},\dots,z_{d})\mathrm{gr}_{\mathfrak{m}}(A/I(S)}\bigg)=\beta_{i}\big(\mathrm{gr}_{\mathfrak{m}}(A/I(S)\big).$$
$A/I(S)$ being Cohen-Macaulay, $z_{1},\dots,z_{d}$ form a regular sequence 
in $A/I(S)$. Hence,
\begin{align*}
\beta_{i}\big(\mathrm{gr}_{\mathfrak{m}}(A/I(S)\big)=&\beta_{i}\big(\bar{A}/\pi_{d}(I(S))\big)\\
=&\beta_{i}\big(A/(z_{1},\dots,z_{d},I(S)\big)\\
=&\beta_{i}\big(A/I(S)\big). \qed
\end{align*} 

\section{Homogeneous Simplicial Affine Semigroup}
The main aim of this section is to generalize the concept of  
homogeneous numerical semigroups to simplicial affine semigroups. 
Let us recall some definitions and examples. Let $S$ be a simplicial 
affine semigroup in $\mathbb{N}^{d}$ minimally generated by 
$a_{1},\dots,a_{d},a_{d+1},\dots,a_{d+r}$, where $a_{1},\dots,a_{d}$ 
are the extremal rays of $S$. 

Given $0 \neq s \in S,$ the set of lengths of $s$ in $S$ is defined as 
\[
\mathcal{T}(s)=\Bigg\{\sum_{i=1}^{d+r} r_{i} \mid s=\sum_{i=1}^{d+r}r_{i}a_{i}, \, r_{i}\geq 0 \Bigg\}.
\]

\begin{definition}{\rm 
A subset $T \subset S$ is called \emph{homogeneous} if either it is empty or 
$\mathcal{T}(s)$ is singleton for all $0 \neq s \in T$. A simplicial affine 
semigroup $S$, with the set of extremal rays 
$E$, is called \textit{homogeneous} if the Ap\'{e}ry set 
$\mathrm{Ap}(S,E)$ is homogeneous. Hence, all the expressions 
of elements of $\mathrm{Ap}(S,E)$ are maximal (see definition \ref{maximal}).
}
\end{definition}

\begin{example}{\rm 
Let $S$ be a simplicial affine semigroup, with extremal rays 
$E=\{a_{1},\dots,a_{d}\}$, such that the defining ideal $I(S)$ 
is generic, i.e., all the variables belong to the support of 
these binomials in $I(S)$. 

We show that every simplicial affine semigroup $S$, with a generic $I(S)$, is homogeneous.  
If $b \in \mathrm{AP}(S,E)$ has two expressions, 
i.e., 
$b=\sum_{j=d+1}^{d+r }p_{j}a_{j}=\sum_{j=d+1}^{d+r }q_{j}a_{j}$, 
with $p_{j}\neq q_{j}$ for some $j$, then 
$0\neq z^{p}-z^{q} \in I(S)$. However, $z_{i}$ does not divide any term of 
$z^{p}-z^{q}$, which is a contradiction as $I(S)$ is generic. 
Hence, every element of $\mathrm{AP}(S,E)$ has a unique expression, therefore $S$ is homogeneous.
}
\end{example}

\begin{lemma}\label{ord} 
The following statements are equivalent:
\begin{enumerate}[(a)]
\item $\mathrm{gr}_{\mathfrak{m}}(k[S])$ is Cohen-Macaulay and $(x^{a_{1}},\dots,x^{a_{d}})$ 
is a reduction ideal of $\mathfrak{m}$;

\item $k[S]$ is Cohen-Macaulay and $\mathrm{ord}_{S}(b + a_{i})= \mathrm{ord}_{S}(b)+1, \,\text{for all}\, b \in S$ and $i = 1,\dots , d$;

\item $k[S]$ is Cohen-Macaulay and $\mathrm{ord}_{S}(b + \sum\limits_{i=1}^{d}a_{i})= \mathrm{ord}_{S}(b)+\sum\limits_{i=1}^{d}n_{i}, \,\text{for all}\, b \in \mathrm{Ap}(S,E)$ and $n_{1},\dots,n_{d}\in \mathbb{N}$.
\end{enumerate}
\end{lemma}

\proof See Proposition 5.4 in \cite{Reduction}. \qed
\medskip

\begin{remark}\label{tors} Let 
$T(S):=\{b \in S \mid \mathrm{ord}_{S}(b+ a_{i})> \mathrm{ord}_{S}(b)+1, \, i\in \{1,\dots,d\}\}$. 
We have $\mathrm{ord}_{S}(b+ a_{i})\geq \mathrm{ord}_{S}(b)+1$,  
therefore by Lemma \ref{ord}, $T(S)=\phi$ if and only if $\mathrm{gr}_{\mathfrak{m}}(R)$ 
is Cohen-Macaulay.
\end{remark}

\noindent\textbf{Notations.} For a tuple $p=(p_{1},\dots,p_{i},\dots,p_{d+r})$, 
we define 
\begin{itemize}
\item $|p|=\sum\limits_{j=1}^{d+r}p_{j}$,

\item $r(p)=\sum\limits_{j=1}^{d+r}p_{j}a_{j}$,

\item For $i \in \{1,\dots,d+r\}$, 
\[  \bar{p} =
    \begin{cases}
      p & \text{if $p_{i}=0$}\\
      (p_{1},\dots,p_{i}-1,\dots,p_{d+r}) & \text{if $p_{i}>0$}
    \end{cases}       
 \]

\end{itemize}
\medskip

\noindent The next Theorem is a generalization of Theorem 3.12 of 
\cite{Homo type}, which was proved in the context of numerical semigrops. 
We show that similar results can be proved for affine simplicial semigroups 
as well.  
We borrow the main ideas from their proof, with the exception that, 
we define the maps $\pi_{d}$ which retain the homogeneity of the 
homogeneous part of $I(S)$ and map the non-homogeneous elements to 
monomials.  

\begin{theorem}\label{homo}
Let $S$ be a simplicial affine semigroup. The following statements 
are equivalent. 
\begin{enumerate}[(a)]
\item  $S$ is homogeneous and $\mathrm{gr}_{\mathfrak{m}}(k[S])$ is Cohen-Macaulay.

\item For all $z^{p}-z^{q} \in I(S)$, with $|p| > |q|$, we have 
$r(p) \notin \mathrm{Ap}(S,E)$. Moreover, if $ \bar{q} $ is a maximal factorization, 
then $p_{i}\geq q_{i}$, for all $i=1,\dots, d$.

\item There exists a minimal generating set of binomials generators $J$ for $I(S)$, 
such that if $z^{p}-z^{q}\in J$ with $|p|>|q|$, then $p_{i} \neq 0$ for some $i=1,\dots,d$.

\item There exists a minimal generating set of binomials generators $J$ for $I(S)$, 
which is a standard basis, and for all $z^{p}-z^{q}\in J$, with $|p|>|q|$,  we have 
$p_{i} \neq 0$ for some $i=1,\dots,d$.

\item There exists a minimal Gr\"{o}bner basis $G$ of $I(S)$, with respect to the 
negative degree reverse lexicographic ordering induced by 
$z_{d+r} > \dots > z_{d} > \dots >z_{1}$, such that for every $i=1,\dots,d$, 
the variable $z_{i}$ does not divide the leading monomial of any element of $G$, 
and there exists $1\leq j \leq d$ such that $z_{j}$ belongs to the support 
of all non-homogeneous elements of $G$.

\end{enumerate}
\end{theorem}

\proof (a) $\Rightarrow$ (b). Let $z^{p}-z^{q} \in I(S)$, with $|p| > |q|$. 
Therefore, $r(p)=\sum\limits_{j=1}^{d+r}p_{j}a_{j}=\sum\limits_{j=1}^{d+r}q_{j}a_{j}=r(q)$. 
By the definition of homogeneous affine semigroups, all the expressions of elements 
of $\mathrm{Ap}(S,E)$ are maximal, however $|p| > |q|$, therefore $r(p) \notin \mathrm{Ap}(S,E)$. 
Let  $s:= r(\bar{q})$, then $\mathrm{ord}_{S}(s)=\sum\limits_{j=1}^{d+r}\bar{q}$. 
Suppose there exists some  $1 \leq i \leq d$, such that $p_{i}<q_{i}$. Then, 
$\bar{q}=(q_{1},\dots,q_{i-1},q_{i}-1,q_{i+1},\dots,q_{d+r})$ 
and $s+a_{i}=\sum\limits_{j=1}^{d+r}q_{j}a_{j}=\sum\limits_{j=1}^{d+r}p_{j}a_{j}$. 
Therefore, 
$$\mathrm{ord}_{S}(s+a_{i})\geq \sum\limits_{j=1}^{d+r}p_{j}>\sum\limits_{j=1}^{d+r}q_{j}=\sum\limits_{j=1,j \neq i}^{d+r}q_{j}+(q_{i}-1)+1= \mathrm{ord}_{S}(s)+1,$$ 
which is a contradiction of Remark \ref{tors}. Hence $p_{i}\geq q_{i}$, for all $i=1,\dots,d$.
\medskip

\noindent (b) $\Rightarrow$ (c). Let $\mathcal{T}_{1}$ be a set of generators for $I(S)$ 
and let $f_{1}=z^{p}-z^{q}\in \mathcal{T}_{1}$. with $|p|>|q|$ and $p_{i}=0$ for all 
$i=1,\dots,d$. Then, 
$s=\sum\limits_{i=1}^{d+r}q_{i}a_{i}=\sum\limits_{i=1}^{d+r}p_{i}a_{i} \notin \mathrm{Ap}(S,E)=\cap_{i=1}^{d}\mathrm{AP}(S,a_{i})$, which implies that 
$s \notin \mathrm{Ap}(S,a_{i})$, for some $i=1,\dots,d$, and 
$s=\sum\limits_{j=1}^{d+r}r_{j}a_{j}$, such that $r_{i}>0$ for some 
$i=1,\dots,d$. Now, $\mathcal{T}_{2}= (\mathcal{T} \setminus \{f_{1}\} \cup \{z^{p}-z^{r},z^{q}-z^{r}\})$ is again a finite set of generators for $I(S)$, such that $r_{i} \neq 0$ for some $i=1,\dots,d$. 
By continuing this way, we get the generating set $\mathcal{T}$ for $I(S)$, 
such that $z^{r}-z^{r'}\in J$ with $|r|>|r'|$  and $r_{i} \neq 0$ for some 
$i=1,\dots,d$. Now a minimal generating set $J$ for $I(S)$, extracted from 
$\mathcal{T}$, has the desired property.
\medskip

\noindent (c) $\Rightarrow$ (d). Let $J=\{f_{1},\dots,f_{r},g_{1},\dots,g_{s}\}$ be 
a minimal generating set of binomials for $I(S)$, where $f_{1},\dots,f_{r}$ are homogeneous 
and $g_{1},\dots,g_{s}$ are non-homogeneous. We consider the map 
$$\pi_{d}:k[z_{1},\dots,z_{d},\dots,z_{d+r}] \rightarrow k[z_{d+1},\dots,z_{d+r}],$$ 
such that $\pi_{d}(z_{j})=0, 1 \leq j \leq d$ and $\pi_{d}(z_{j})=z_{j}, d+1 \leq j \leq d+r$. 

Let $g=z^{p'}-z^{q'}$, with $|p'|>|q'|$. Then from (c), $p'_{i}\neq 0$ for some 
$i=1,\dots, d$, and $z_{i} \mid z^{p'}$ implies that $\pi_{d}(g)=z^{q'}$. 
Therefore, 
$$E=\{\pi_{d}(f_{1}),\dots,\pi_{d}(f_{r}),\pi_{d}(g_{1}),\dots,\pi_{d}(g_{s})\}$$ 
generates $\pi_{d}(I(S))$. Since $E$ is a set of homogeneous set of generators of 
$\pi_{d}(I(S))$, it is a standard basis of $\pi_{d}(I(S)).$ 
From (\cite{Super}, Theorem 1), $J$ is a standard basis of $I(S)$.
\medskip

\noindent (d) $\Rightarrow$ (e). Follows from  Theorem 3.12 in \cite{Homo type}.
\medskip

\noindent (e) $\Rightarrow$ (a). Suppose $b \in \mathrm{AP}(S,E)=\cap_{i=1}^{d}\mathrm{AP}(S,a_{i})$, such that $b=\sum\limits_{i=1}^{d+r }p_{i}a_{i}=\sum\limits_{i=1}^{d+r }q_{i}a_{i}$. This implies $z^{p}-z^{q}$ is homogeneous, otherwise, from the hypothesis $z_{i}$ must belong to the support of $z^{p}-z^{q}$, 
for some $ 1 \leq i \leq d$. Assume that $z_{i}$ divides $z^{p}$, 
then $z^{p}=z_{1}^{p_{1}}\dots z_{i}^{p_{i}}\dots z_{d+r}^{p_{d+r}}$, 
with $p_{i} \geq 1$. We have
$b=\sum\limits_{i=1}^{d+r }p_{i}a_{i}=p_{1}a_{1}+\dots+(p_{i}-1)a_{i}+a_{i}+\dots+p_{d+r}a_{d+r}$. Hence 
$b-a_{i}=p_{1}a_{1}+\dots+(p_{i}-1)a_{i}+\dots+p_{d+r}a_{d+r}\in S$, therefore 
$b \notin \mathrm{AP}(S,a_{i})$, for all $ 1 \leq i \leq d$, which is a 
contradiction as $b \in \mathrm{AP}(S,E)$. Therefore, 
$|p|=|q|$ and $\mathrm{AP}(S,E)$ is homogeneous and from Theorem \ref{criterion}. 
Hence, $\mathrm{gr}_{\mathfrak{m}}(k[S])$ is Cohen-Macaulay. \qed

\begin{definition}[\cite{Homo type}, Definition 3.14]
\emph{A semigroup $S$ is said to be \emph{of homogeneous type} if 
$\beta_{i}(k[S])=\beta_{i}(\mathrm{gr}_{\mathfrak{m}}(k[S]))$ for all $i \geq 1$}. 
\end{definition}

\begin{theorem}\label{Betti}
Let $S$ be a simplicial affine homogeneous semigroup such that $\mathrm{gr}_{\mathfrak{m}}(k[S])$ is Cohen-Macaulay. Then $\beta_{j}(k[S])=\beta_{j}(\mathrm{gr}_{\mathfrak{m}}(k[S]))$, for all $j \geq 1$.
 
\end{theorem}

\proof $S$ is a simplicial affine homogeneous semigroup such that $\mathrm{gr}_{\mathfrak{m}}(k[S])$ is Cohen-Macaulay. Therefore, by Theorem \ref{homo}, there exists a minimal Gr\"{o}bner basis $G$ of $I(S)$ with 
respect to the negative degree reverse lexicographic ordering induced by 
$z_{d+r} > \dots > z_{d} > \dots >z_{1}$, with the following properties: 
$z_{j}$ does not divide the leading monomial 
of any element of $G$, for every $1 \leq j \leq d$, and there exists $1\leq j \leq d$, 
such that $z_{j}$ 
belongs to the support of all non-homogeneous elements of $G$. Hence, by Theorem \ref{CM}, we 
can write $\beta_{j}(k[S])=\beta_{j}(\mathrm{gr}_{\mathfrak{m}}(k[S]))$, for all $ j \geq 1$.  \qed

\begin{remark}{\rm
If $S$ is a simplicial affine semigroup of homogeneous type, such that $k[S]$ is Cohen-Macaulay, 
then $\mathrm{depth}(\mathrm{gr}_{\mathfrak{m}}(k[S]))=\mathrm{depth}(k[S]) = 
\mathrm{dim}(k[S])=\mathrm{dim}(\mathrm{gr}_{\mathfrak{m}}(k[S]))$ (see exercise 13.8, \cite{Eisenbud}). 
Hence, $\mathrm{gr}_{\mathfrak{m}}(k[S])$ is Cohen-Macaulay.
}
\end{remark}

\begin{example}{\rm 
(Example 4.12, \cite{Reduction}) 
Assume that $S$ is generated by $a_{1}=(0,2),a_{2}=(2,1),a_{3}=(0,3)$, and $a_{4}=(1,2)$, 
with extremal rays $a_{1},a_{2}$. Then  $\mathrm{gr}_{\mathfrak{m}}(k[S])$ is Cohen-Macaulay and 
$\mathrm{AP}(S,E)=\{(0,0),(0,3),(1,2),(1,5)\}$. Note that every element of 
$\mathrm{AP}(S,E)$ has a unique expression, hence $S$ is of homogeneous type.
}
\end{example}

\begin{example}{\rm 
Backelin defined the class of semigroups $\langle s, s+3, s+3n+1, s+3n+2\rangle$, 
for $n\geq 2$, $r\geq 3n+2$ and $s=r(3n+2)+3$. Let 
$\tilde{S}=\langle (0,s+3n+2),(s,3n+2),(s+3,3n-1),(s+3n+1,1),(s+3n+2,0)\rangle \subset \mathbb{N}^{2}$. 
It is known that $k[\tilde{S}]$ is Cohen-Macaulay (see \cite{Backelin}, Theorem 2.9). 
Note that $\{(0,s+3n+2),(s+3n+2,0)\}$ is the set of extremal rays of $\tilde{S}$ and 
$z_{0}, z_{4}$ belong to the support of non-homogeneous elements of a Gr\"{o}bner 
basis of the defining ideal of the projective closure of Backelin's curve (see \cite{Backelin}, Theorem 2.5). 
Hence $\tilde{S}$ is of homogeneous type. 
}
\end{example}

\section{Nice Extension Of Simplicial Affine Semigroup}
In this section, we develop the concept of the nice extension of simplicial affine 
semigroups, which is a generalization of the nice extension of numerical semigroups 
given in \cite{Nice}.

\begin{definition} {\rm 
Let $S$ be a simplicial affine semigroup, fully embedded in $\mathbb{N}^{d}$, 
minimally generated by $a_{1},\dots,a_{d},a_{d+1},\dots,a_{d+r}$, such that 
$a_{1},\dots,a_{d}$ are the extremal rays of $S$. Suppose $b \in \langle S \rangle $ 
and  $\lambda, \mu \in \mathbb{N}$, with $\mathrm{gcd}(\lambda, \mu)=$1. The 
semigroup $S_{b}=\lambda S \cup  \{\mu b\}$ is an extension of $S$. Let 
$b =\alpha_{1}a_{1}+\alpha_{2}a_{2}+\cdots +\alpha_{d+r}a_{d+r}$, where 
$\alpha_{1}, \dots, \alpha_{d} \in \mathbb{N}$. If $\lambda \leq \sum\limits_{i=1}^{d+r}\alpha_{i}$, 
then $S_{b}$ is called the \emph{nice extension} of $S$.
}
\end{definition}

 \begin{remark}
We write $b =\alpha_{1}a_{1}+\alpha_{2}a_{2}+\cdots +\alpha_{d+r}a_{d+r}$, where 
$\alpha_{1}, \dots, \alpha_{d} \in \mathbb{N}$. By Proposition 1 in \cite{GB}, 
the defining ideal of $k[S_{b}]$ is $I(S_{b})=I(S) \cup \{y^{\lambda}-z_{1}^{\mu \alpha_{1}}\dots z_{d+r}^{\mu \alpha_{d+r}}\}$. Therefore, $\mu(I(S_{b}))=1+\mu(I(S))$, where $\mu(I(S_{b}))$ 
and $\mu(I(S))$ denote 
the minimal number of generators of the ideals $I(S_{b})$ and $I(S)$ respectively. 
\end{remark}

\begin{lemma}\label{extremal}
Let $S$ be a simplicial affine semigroup, fully embedded in $\mathbb{N}^{d}$, minimally  generated 
by $a_{1},\dots,a_{d}
,a_{d+1},\dots,a_{d+r}$, such that $a_{1},\dots,a_{d}$ are the extremal rays of $S$. Then 
the extension $S_{b}$ is a simplicial affine semigroup minimally generated by $\lambda a_{1},\dots,\lambda a_{d}
,\lambda a_{d+1},\dots,\lambda a_{d+r},\mu b$, with extremal rays $\{\lambda a_{1},\dots,\lambda a_{d}\}$. 
\end{lemma}

\proof Since $b \in S$, there exist $ q_{1},\dots,q_{d} \in \mathbb{Q}$, such that 
$b = q_{1}a_{1}+\dots +q_{d}a_{d}$. Therefore, 
$\mu b = \frac{\mu q_{1}}{\lambda}(\lambda a_{1})+\dots+\frac{\mu q_{d}}{\lambda}(\lambda a_{d})$. 
Hence, $\{\lambda a_{1},\dots,\lambda a_{d}\}$ is the set of extremal rays of $S_{b}$ and $S_{b}$ 
is a simplicial affine semigroup. \qed

\begin{theorem}\label{Complete}
Let $S_{b}$ be an extension of a simplicial affine semigroup $S$, with affine semigroup rings 
$k[S_{b}]$ and $k[S]$ respectively. If $k[S]$ is a complete intersection then $k[S_{b}]$ is 
also a complete intersection.
\end{theorem}

\proof Since $ b \in \langle S \rangle$, $\mathrm{cone}(S_{b})$ generates the same subspace as $\mathrm{cone}(S)$, therefore $\mathrm{dim}(k[S_{b}])=\mathrm{dim}(k[S])$. 
Now, since $k[z_{1},\dots,z_{d+r},y]$ is a regular ring, we have
\begin{align*}
\mathrm{ht}(I(S_{b}))&=\mathrm{dim}(k[z_{1},\dots,z_{d+r},y])-\mathrm{dim}\big(\frac{k[z_{1},\dots,z_{d+r}]}{I(S_{b})}\big)\\
&=(d+r)+1-\mathrm{dim}(k[S_{b}])\\
&=(d+r)+1-\mathrm{dim}(k[S])\\
&=(d+r)-\mathrm{dim}(k[S])+1\\
&=\mathrm{ht}(S)+1\\
&=\mu (I(S))+1 \,\,\,
\\&=\mu (I(S_{b}).
\end{align*}
Therefore, $k[S_{b}]$ is also a complete intersection. \qed

\begin{definition}{\rm 
A simplicial affine semigroup $S$ in $\mathbb{N}^{d}$ is obtained by a sequence of nice extensions 
if there are affine semigroup $S^{0},\dots,S^{l}$, such that $S^{0}$ is the semigroup generated by 
$\{(1,0,\dots,0),(0,\dots,0,1)\}$, $S^{l}=S$ and $S^{i+1}$ is a nice extension of $S^{i}$, for every 
$i=0,\dots,l-1$ and for some $l \in \mathbb{N}$. 
}
\end{definition}

\begin{theorem}\label{CI}
Every semigroup ring associated with an affine semigroup, obtained by a sequence of 
nice extensions, is a complete intersection.
\end{theorem}

\proof The proof is by induction. For $i=0$, the semigroup ring $k[\mathbb{N}^{d}]$ 
is isomorphic to a polynomial ring, therefore $k[\mathbb{N}^{d}]$ is a complete intersection. 
Let the statement be true for $i=r$, i.e, let $k[S^{r}]$ be a complete intersection. 
Since $S^{r+1}$ is an extension of $S^{i}$ and $k[S^{r}]$ is a complete intersection by 
induction hypothesis, it follows from Theorem \ref{Complete} that $k[S^{r+1}]$ is a 
complete intersection. \qed

\begin{theorem}\label{CM glu} 
Let $S$ be a simplicial affine semigroup in $\mathbb{N}^{d}$, minimally  generated by $a_{1},\dots,a_{d}
,a_{d+1},\dots,a_{d+r}$, such that $a_{1},\dots,a_{d}$ are the extremal rays of $S$. 
Let us assume 
that the associated graded ring $\mathrm{gr}_{\mathfrak{m}}(k[S])$ is Cohen-Macaulay. 
Let $S_{b}$ be a nice extension of $S$, then the associated graded ring $\mathrm{gr}_{\mathfrak{m}_{b}}(k[S_{b}])$ 
is Cohen-Macaulay. 
\end{theorem}

\proof Let $G=\{f_{1},\cdots,f_{r}\}$ be a minimal Gr\"{o}bner basis of the defining ideal $I(S)$ 
of the semigroup ring $k[S]$, with respect to the negative degree reverse lexicographic 
ordering induced by $z_{d+r} > \dots > z_{d} > \dots >z_{1}$. We claim that 
$G_{b}=\{f_{1},\dots,f_{r},y^{\lambda}-z_{1}^{\mu \alpha_{1}}\dots z_{d+r}^{\mu \alpha_{d+r}}\}$ is a minimal Gr\"{o}bner basis of the defining ideal $I(S_{b})=I(S) \cup \{y^{\lambda}-z_{1}^{\mu \alpha_{1}}\dots z_{d+r}^{\mu \alpha_{d+r}}\}$ of the semigroup ring $k[S_{b}]$, 
with respect to the monomial order written above. Since $S_{b}$ is a nice 
extension of $S$ and $\lambda \leq \sum\limits_{i=1}^{d+r}\alpha_{i}$, therefore $\mathrm{LM}(y^{\lambda}-z_{1}^{\mu \alpha_{1}}\dots z_{d+r}^{\mu \alpha_{d+r}})=y^{\lambda}$. 
We note that $y$ does not appear in any $f_{i}$, for $1 \leq i \leq r$, 
and the leading monomials $\mathrm{LM}(f_{i})$ and 
$\mathrm{LM}(y^{\lambda}-z_{1}^{\mu \alpha_{1}}\dots z_{d+r}^{\mu \alpha_{d+r}})$ 
are mutually coprime, therefore, the $S$-polynomial 
$S(f_{i},y^{\lambda}-z_{1}^{\mu \alpha_{1}}\dots z_{d+r}^{\mu \alpha_{d+r}})$ reduces to zero
when divided by $G_{b}$. Also, $G$ is a minimal Gr\"{o}bner basis, therefore 
$S(f_{i},f_{j})$ reduces to zero upon division by $G$ and hence upon division by 
$G_{b}$. By the 
Buchberger's criterion, the set  
$G_{b}=\{f_{1},\cdots,f_{r},y^{\lambda}-z_{1}^{\mu \alpha_{1}}\dots z_{d+r}^{\mu \alpha_{d+r}}\}$ is a minimal Gr\"{o}bner basis of the defining ideal $I(S_{b})$ of the semigroup ring $k[S_{b}]$, 
with respect to the said order. From Lemma \ref{extremal}, 
$\lambda a_{1},\dots,\lambda a_{d}$ are also the extremal rays of $S_{b}$ and 
since $\mathrm{gr}_{\mathfrak{m}}(k[S])$ is Cohen-Macaulay, it follows from Theorem \ref{criterion} 
that for every $j=1,\dots,d$, the indeterminate $z_{j}$ does not divide 
any element of $G$ and it does not divide 
$(y^{\lambda}-z_{1}^{\mu \alpha_{1}}\dots z_{d+r}^{\mu \alpha_{d+r}})$. 
Therefore, 
for every $1 \leq j \leq d$, the indeterminate $z_{j}$ does not divide any element 
of $G_{b}$, hence $\mathrm{gr}_{\mathfrak{m}_{b}}(k[S_{b}])$ is Cohen-Macaulay by 
Theorem \ref{criterion}. \qed

\begin{corollary}
Let $S$ be a homogeneous simplicial affine semigroup in $\mathbb{N}^{d}$, minimally 
generated by $a_{1},\dots,a_{d}
,a_{d+1},\dots,a_{d+r}$, such that $a_{1},\dots,a_{d}$ are the extremal rays of $S$. 
Then, the nice extension $S_{b}$ of $S$, for $b \in \langle S \rangle$, 
is also a homogeneous simplicial affine semigroup.
\end{corollary}

\proof From the proof of Theorem \ref{CM glu}, it is clear that $I(S_{b})$ has a 
minimal Gr\"{o}bner basis with respect to the negative degree reverse lexicographic 
ordering induced by $z_{d+r} > \dots > z_{d} > \dots >z_{1}$. Moreover, 
$z_{j}$ does not divide the leading monomial of any element of $G$, for 
every $j=1,\dots,d$, and there exists $j$, $1\leq j \leq d$, 
such that $z_{j}$ belongs to the support of all non-homogeneous elements of $G$. 
Hence, $S_{b}$ is homogeneous by Theorem \ref{homo}(e).\qed

\begin{theorem}
Let $S$ be a homogeneous simplicial affine semigroup in $\mathbb{N}^{d}$, minimally 
generated by $a_{1},\dots,a_{d}, a_{d+1},\dots, a_{d+r}$, such that $a_{1},\dots,a_{d}$ 
are the extremal rays of $S$. Let $S_{b}$ be a nice extension of $S$, for 
$b \in \langle S \rangle$. Then for all $i \geq 1$
\[
\beta_{i}(k[S_{b}])=\beta_{i}(k[S])+\beta_{i-1}(k[S]).
\]
\end{theorem}

\proof Follows from Theorem 1 in \cite{GB}.\qed
 
\section{Numerical Semigroup minimally generated by Geometric sequence}
In this section, we present a particular class of numerical semigroup and its projective 
closure (an affine semigroup) as an illustration of some of the theorems proved in the 
earlier sections. 
\medskip

Let $\mathrm{gcd}(a,b)=1, a <b$ and $r \in \mathbb{N}$. Consider a numerical semigroup $S$ minimally generated by $m_{1}=a^{r}<m_{2}=a^{r-1}b<\cdots<m_{r}=ab^{r-1}<m_{r+1}=b^{r}$. 
Let $k$ be a field and $k[S]:=k[t^{m_{1}},\dots,t^{m_{r+1}}]\subset k[t]$ be the numerical semigroup ring defined by $S$. Let $\eta:A=k[z_{1},\dots,z_{r+1}] \rightarrow k[t]$ be 
the mapping defined by $\eta(z_{i})=t^{m_{i}}, 1 \leq i \leq r+1$. Then, 
$\frac{A}{\mathrm{Ker} (\eta)}\cong k[S]$ is the coordinate ring of the 
affine monomial curve in $\mathbb{A}_{k}^{r+1}$ and $\mathrm{Ker} (\eta)$ is the 
defining ideal of that curve denoted by $\mathfrak{p}$. Let $\mu(\mathfrak{p})$ 
denotes the minimal number of generator of $\mathfrak{p}$.

\begin{theorem}\textbf{(Gastinger)} \label{gastinger}\cite{Gastinger}
Let $A = k[z_{1},\ldots,z_{r}]$ be the polynomial ring, $I\subset A$ 
the defining ideal of a monomial curve defined by natural numbers 
$a_{1},\ldots,a_{r}$, whose greatest common divisor is $1$.  
Let $J$ be an ideal contained in $I$. Then $J = I$ if and 
only if $\mathrm{dim}_{k} A/\langle J + (z_{i}) \rangle =a_{i}$, 
for some $i$; equivalently 
$\mathrm{dim}_{k} A/\langle J + (z_{i}) \rangle =a_{i}$ for any $i$.
\end{theorem}  

\begin{theorem}
The defining ideal $\mathfrak{p}$ of 
the monomial curve defined by $S$, with the coordinate ring $k[S]$, 
is minimally generated by following 
set of binomials 
$$\{P_{1}=z_{2}^{a}-z_{1}^{b},P_{2}=z_{3}^{a}-z_{2}^{b},\dots,P_{r}=z_{r+1}^{a}-z_{r}^{b}\}.$$
\end{theorem}

\proof Let $g_{i}=z_{i}^{b-a}g_{i-1}+P_{i}=z_{i+1}^{a}-z_{1}^{b}z_{2}^{b-a}\dots z_{i}^{b}$, 
for $1 \leq i \leq r$.
Consider $I=\langle g_{1},\dots,g_{r} \rangle$. Then, 
$I \subset \mathfrak{p}$ and 
$A/\langle I+(z_{1})\rangle=\langle z_{2}^{a},z_{3}^{a},\dots,z_{r+1}^{a} \rangle$ 
(in $k[z_{2},\dots,z_{r+1}]$) is a vector space over $k$ with a  
basis consisting of the images of monomials 
$z_{2}^{i_{1}}z_{3}^{i_{2}}\dots z_{r+1}^{i_{r}}$, 
where $0 \leq i_{1},i_{2},\dots,i_{r} \leq a-1$. Therefore, 
$\mathrm{dim}_{k}A/\langle I+(z_{1})\rangle=a^{r}$. Hence, 
$I=\mathfrak{p}$ and since 
$I \subseteq \langle P_{1},\dots,P_{r} \rangle \subseteq \mathfrak{p}$, it 
follows that $\mathfrak{p}=\langle P_{1},\dots,P_{r} \rangle$. \qed

\begin{theorem}\label{grobner}
Let us consider the negative degree reverse lexicographic monomial order 
on $k[z_{1},\dots,z_{r+1}]$, induced by $z_{r+1}> \cdots > z_{1}$. Then 
$G=\{P_{1},\dots,P_{r}\}$ is a minimal Gr\"{o}bner basis of the defining ideal 
$\mathfrak{p}$ of the monomial curve defined by $S$, with the coordinate ring $k[S]$.
\end{theorem}

\proof The leading monomial of $P_{i}$'s are $\mathrm{LM}(P_{i})=z_{i+1}^{a}$, 
for $i=1,\dots,r$, with respect to the given monomial order and 
$\gcd(P_{i},P_{j})=1$, for $i \neq j$. Therefore, the $S$-polynomial 
$S(P_{i},P_{j})$ reduces to zero upon division by $G$. Hence, 
$G$ is a minimal Gr\"{o}bner basis of the ideal $\mathfrak{p}$. \qed

\begin{corollary}
The associated graded ring $\mathrm{gr}_{\mathfrak{m}}(k[S])$ of $k[S]$ is Cohen-Macaulay.
\end{corollary}
\proof From Theorem \ref{grobner}, note that $z_{1}$ does not divide the leading monomial of any elemenet of $G$. The result follows from Theorem \ref{criterion}. \qed 
\medskip

We now discuss about the projective closure of $k[S]$. Consider a map 
$\eta^{h}:k[z_{0},\dots,z_{r+1}]\rightarrow k[u,v]$, 
such that $\eta^{h}(z_{0})=v^{m_{r+1}},\eta^{h}(z_{i})=u^{m_{i}}v^{m_{r+1}-m_{i}}$, 
for $i=1,\dots,r+1$. Then, the homogenization of the ideal $\mathfrak{p}$, 
with respect to the variable $z_{0}$ is $\overline{\mathfrak{p}}$. Thus, 
the projective curve 
$\{[(a^{m_{r+1}}:a^{m_{r+1}-m_{1}}b^{m_{1}}:\cdots:b^{m_{r+1}})]\in\mathbb{P}^{r+1}_{k}\mid a,b\in k\}$ 
is the projective closure of the affine curve 
$k[S]=\{(b^{m_{1}},\ldots b^{m_{r+1}})\in \mathbb{A}^{r+1}_{k}\mid b\in k\}$, 
denoted by $\overline{k[S]}$. 

\begin{theorem}
The rings $k[S]$ and its projective closure $\overline{k[S]}$ 
are both complete intersections.
\end{theorem}

\proof The height of the defining ideal $\mathfrak{p}$ is 
$$\mathrm{ht}(\mathfrak{p}) = 
\mathrm{dim}(k[z_{1},\dots,z_{r+1}])-\mathrm{dim}\big(\frac{k[z_{1},\dots,z_{r+1}]}{\mathfrak{p}}\big) = 
r = \mu(\mathfrak{p}).$$
Therefore $k[S]$ is a complete intersection. Similarly, it can be proved that 
the projective closure $\overline{k[S]}$ of the monomial curve $k[S]$ is also complete intersection.\qed

\begin{corollary}
The projective closure $\overline{k[S]}$ of $k[S]$ is Cohen-Macaulay and Gorenstein.
\end{corollary}

It is a well-known theorem in Commutative Algebra (see Theorem 21.2 in \cite{M}) 
that a local, Noetherian ring is a complete intersection if and only if it can 
be written as a quotient of a regular local ring by a regular sequence. 
It follows from the above observations that 
$\{ P_{1},\dots,P_{r} \}$ form a $A$-regular sequence and therefore 
the defining ideal 
$\mathfrak{p}$ is minimially resolved by the Koszul complex. The 
Betti numbers are give by $\beta_{i}^{A}(A/\mathfrak{p})={r \choose i}$.

\bibliographystyle{amsalpha}

\end{document}